\documentclass[12pt]{article}
\usepackage{amsmath,amssymb,amsfonts,amsthm}
\usepackage[all]{xy}
\newcounter{item}[section]
\newcounter{kirshr}
\newcounter{kirsha}
\newcounter{kirshb}
\newenvironment{enumroman}{\setcounter{kirshr}{1}
\begin{list}{(\roman{kirshr})}{\usecounter{kirshr}} }{\end{list}}
\newenvironment{enumarab}{\setcounter{kirshb}{1}
\begin{list}{(\arabic{kirshb})}{\usecounter{kirshb}} }{\end{list}}
\newtheorem{theorem}{Theorem}[section]

\newtheorem{lemma}[theorem]{Lemma}
\newtheorem{corollary}[theorem]{Corollary}
\newenvironment{demo}[1]{\noindent{\bf #1.}\upshape\mdseries}
{\nopagebreak{\hfill\rule{2mm}{2mm}\nopagebreak}\par\normalfont}
\theoremstyle{definition}

\newtheorem{example}[theorem]{Example}
\newtheorem{definition}[theorem]{Definition}

\def\C{{\mathfrak{C}}}

\def\Nr{{\mathfrak{Nr}}}
\def\Fr{{\mathfrak{Fr}}}
\def\Sg{{\mathfrak{Sg}}}

\def\A{{\mathfrak{A}}}
\def\B{{\mathfrak{B}}}
\def\C{{\mathfrak{C}}}
\def\D{{\mathfrak{D}}}

\def\CA{{\bf CA}}

\def\FPA{{\bf FPA}}
\def\Lf{{\bf Lf}}
\def\Dc{{\bf Dc}}

\def\RK{{\bf RK}}
\def\RCA{{\bf RCA}}

\def\Rd{{\mathfrak{Rd}}}
\def\(R)RA{{\bf (R)RA}}

\def\Dc{{\bf Dc}}

\def\Dc{{\bf Dc}}

 \def\CA{{\sf CA}}
\def\B{{\sf B}}
\def\G{{\sf G}}

\def\Nr{{\mathfrak{Nr}}}

\def\Nr{{\mathfrak{Nr}}}

\def\A{{\mathfrak{A}}}
\def\B{{\mathfrak{B}}}
\def\C{{\mathfrak{C}}}
\def\D{{\mathfrak{D}}}

\def\CA{{\bf CA}}

\def\RCA{{\bf RCA}}
\def\G{{\bf G}}

\def\L{{\mathfrak{L}}}

\def\PA{{\bf PA}}

\def\Nr{{\mathfrak{Nr}}}
\def\Fr{{\mathfrak{Fr}}}
\def\Sg{{\mathfrak{Sg}}}

\def\Rd{{\mathfrak{Rd}}}
\def\Ig{{\mathfrak{Ig}}}
\def\CA{{\bf CA}}
\def\RCA{{\bf RCA}}

\def\L{{\bf L}}

\def\(R)RA{{\bf (R)RA}}

\def\Dc{{\bf Dc}}

\def\Dc{{\bf Dc}}

\def\Nr{{\mathfrak{Nr}}}
\def\Fr{{\mathfrak{Fr}}}
\def\Sg{{\mathfrak{Sg}}}

\def\Rd{{\mathfrak{Rd}}}
\def\Ig{{\mathfrak{Ig}}}
\def\CA{{\bf CA}}
\def\RCA{{\bf RCA}}

\def\L{{\bf L}}

\def\RK{{\bf RK}}
\def\(R)RA{{\bf (R)RA}}

\def\Dc{{\bf Dc}}

\def\Dc{{\bf Dc}}

\def\PA{{\bf PA}}

\def\Dc{{\bf Dc}}

 \def\CA{{\sf CA}}

\def\Dc{{\bf Dc}}

\def\G{{\mathfrak{G}}}

\def\A{{\mathfrak{A}}}
\def\B{{\mathfrak{B}}}
\def\C{{\mathfrak{C}}}
\def\D{{\mathfrak{D}}}

\def\Nr{{\mathfrak{Nr}}}
\def\F{{\mathfrak{F}}}
\def\CA{{\bf CA}}
\def\RCA{{\bf RCA}}

\def\Sg{{\mathfrak Sg}}

\title{Neat embeddings as adjoint situations}
\author{Tarek Sayed Ahmed}

\begin{document}
\maketitle

\begin{abstract} Looking at the operation of forming neat $\alpha$-reducts as a functor, with $\alpha$ an infinite ordinal,
we investigate when such a functor
obtained by truncating $\omega$ dimensions, has a right adjoint. 
We show that the neat reduct functor for representable cylindric algebras does not have a right adjoint, while that of polyadic algebras 
is an equivalence. We relate this categorical result to several amalgamation properties for classes of representable algebras.
\end{abstract}
\section{Introduction}
We follow the notation of \cite{HMT1} often without warning. All ordinals considered are infinite.
$\CA_{\alpha}$ denotes the class of cylindric algebras of dimension $\alpha$.
The subclasses $\Lf_{\alpha}$ and $\Dc_{\alpha}$ of locally finite dimensional and dimension complemented algebras, of dimension $\alpha$, respectively,
are defined in \cite{HMT1}, definition 1.11.1.
The class $\Dc_{\alpha}$ is a non-trivial generalization of the class $\Lf_{\alpha}$ and it shares some of its nice properties.
For example if  $\A\subseteq \Nr_{\alpha}\B$ and $A$ generates $\B$ then $\A=\Nr_{\alpha}\B$
and if $\A\subseteq \Nr_{\alpha}\B_1$ and $\A\subseteq \Nr_{\alpha}\B_2$, and $A$ generates both $\B_1$ and $\B_2$
then $\B_1\cong \B_2$ with an isomorphism that is equal to the identity on $\A,$ cf. \cite{HMT1} 2.6.67, 2.6.72.
In \cite{HMT2}, it was posed as a question by Henkin and Monk whether such nice results generalize to the class of representable algebras.
In \cite{n} a negative answer is given.
Here we present this result in a categorical setting by showing that the neat reduct operator viewed as functor has no right adjoint.

We give a contrasting result for polyadic algebras.  We prove that in this case the neat reduct functor
is an equivalence.  We note that the literature for polyadic algebras is extensive \cite{ANS}, \cite{AU}, \cite{Bulletin}, \cite{SUP}, \cite{Hung}, 
\cite{D}, \cite{DM}, with interest converging to their reducts
\cite{FM}, \cite{F1}, \cite{F2}, \cite{trans}.  We start by reviewing some categorical concepts. We follow \cite{cat} for categorical, notions, 
definitions and theorems.
In particular, for a category $L,$ $Ob(L)$ denotes the class of objects of the category and $Mor(L)$ denotes the corresponding class of morphisms.

\begin{definition} Let $L$ and $K$ be two categories.
Let $G:K\to L$ be a functor and let $\B\in Ob(L)$. A pair $(u_B, \A_B)$ wth $\A_B\in Ob(K)$ and $u_B:\B\to G(\A_B)$ is called a universal map
with respect to $G$
(or a $G$ universal map) provided that for each $\A'\in Ob(K)$ and each $f:\B\to G(\A')$ there exists a unique $K$ morphism
$\bar{f}: \A_B\to \A'$ such that
$$G(\bar{f})\circ u_B=f.$$
\end{definition}

\begin{displaymath}
    \xymatrix{
        \mathfrak{B} \ar[r]^{u_B} \ar[dr]_f & G(\mathfrak{A}_\mathfrak{B}) \ar[d]^{G(f)}  &\mathfrak{A}_\mathfrak{B} \ar[d]^{\hat{f}} \\
             & G(\mathfrak{A}')  & \mathfrak{A}'}
\end{displaymath}

The above definition is strongly related to the existence of adjoints of functors.
For undefined notions in the coming definition, the reader is referred to \cite{cat}
Theorem 27.3 p. 196.
\begin{theorem} Let $G:K\to L$.
\begin{enumarab}
\item If each $\B\in Ob(K)$ has a $G$ universal map $(\mu_B, \A_B)$, then there exists a unique adjoint situation $(\mu, \epsilon):F\to G$
such that $\mu=(\mu_B)$ and for each $\B\in Ob(L),$
$F(\B)=\A_B$.
\item Conversely, if we have an adjoint situation $(\mu,\epsilon):F\to G$ then for each $\B\in Ob(K)$ $(\mu_B, F(\B))$ have a $G$ universal map.
\end{enumarab}
\end{theorem}

Now we apply this definition to the `neat reduct functor' from a
certain subcategory of $\CA_{\alpha+\omega}$
to $\RCA_{\alpha}$. More precisely, let 
$$\L=\{\A\in \CA_{\alpha+\omega}: \A=\Sg^{\A}\Nr_{\alpha}\A\}.$$
Note that $\L\subseteq \RCA_{\alpha+\omega}$. The reason is that any $\A\in \L$ is generated by $\alpha$ -dimensional elements,
so is dimension complemented (that is $\Delta x\neq \alpha$ for all $x$), and such algebras are representable.
Consider $\Nr_{\alpha}$ as a functor from $\bold L$ to $\CA_{\alpha}$, but we restrict morphisms to one to one homomorphisms; that is we take only
embeddings.
By the neat embedding theorem $\Nr_{\alpha}$ is a functor from $\L$ to $\RCA_{\alpha}$.
(For when $\A\in \CA_{\alpha+\omega},$ then $\Nr_{\alpha}\A\in \RCA_{\alpha}$).
The question we adress is: Can this functor be ``inverted".
This functor is not dense since there are representable algebras not in $\Nr_{\alpha}\CA_{\alpha+\omega}$, 
as the following example, which is a straightforward  adaptation  
of a result in \cite{SL} shows:
\begin{example}\label{ex} 
\begin{enumarab}
\item  Let $\F$ be a field of characteristic $0$. Let 
$$V=\{s\in {}^{\alpha}\F: |\{i\in \alpha: s_i\neq 0\}|<\omega\},$$
Let
$${\C}=(\wp(V),
\cup,\cap,\sim, \emptyset , V, {\sf c}_i,{\sf d}_{ij})_{i,j\in \alpha},$$
with cylindrifiers and diagonal elements restricted to $V$.  
Let $y$ denote the following $\alpha$-ary relation:
$$y=\{s\in V: s_0+1=\sum_{i>0} s_i\}.$$
Note that the sum on the right hand side is a finite one, since only 
finitely many of the $s_i$'s involved 
are non-zero. 
For each $s\in y$, we let 
$y_s$ be the singleton containing $s$, i.e. $y_s=\{s\}.$ 
Define 
${\A}\in \CA_{\alpha}$ 
as follows:
$${\A}=\Sg^{\C}\{y,y_s:s\in y\}.$$
Then it is proved in \cite{SL} that  
$$\A\notin \Nr_{\alpha}\CA_{\alpha+1}.$$ 
That is for no $\mathfrak{P}\in \CA_{\alpha+1}$, it is the case that $\Sg^{\C}\{y,y_s:s\in y\}$ 
exhausts the set of all $\alpha$ dimensional elements 
of $\mathfrak{P}$.  
\item Let $\A$ be as in above. Then since $\A$ is a weak set algebra, it is representable.
Hence $\A\in S\Nr_{\alpha}\CA_{\alpha+\omega}$. Let $\B\in \CA_{\alpha+\omega}$ be an algebra such that $\A\subseteq \Nr_{\alpha}\B$. 
Let $\B'$ be the subalgebra of $\B$ 
generated by $\A$. Then $\A$ generates $\B$ but $\A$ is not isomorphic to $\Nr_{\alpha}\B$.
\end{enumarab}
\end{example}

Item (2) in the above example says that there are two non isomorphic algebras, namely $\A$ and $\Nr_{\alpha}\B'$ 
that generate the same algebra $\B'$ using extra dimensions \cite{n}.
If $\A\subseteq \Nr_{\alpha}\B$ then $\B$ is called a dilation of $\A$. $\B$ is a minimal dilation if $\A$
generates $\B$, in which case $\A$ is called a generating subreduct of $\B$. In the previous example $\A$ is a generating subreduct of $\B$.
One would expect that the ``inverse" of the Functor $\Nr$
would be the functor that takes $\A$ to a minimal dilation, and lifting morphisms.
But this functor is not even a right adjoint.

A concise formulation of the above is:

\begin{theorem}\label{amal} For $\alpha\geq \omega$, the following hold:
\begin{enumroman}
\item There exist $\A, \A' \in \RCA_{\alpha}$, $\B, \B' \in \CA_{\alpha+\omega}$ with
embeddings $e_A:\A \to \Nr_\alpha \B$ and $e_{A'}:\A' \to \Nr_\alpha \B'$
such that $ \Sg^\B e_A(A) = \B$ and $ \Sg^{\B'} e_{A'}(A) = \B'$, and an isomorphism $ i : \A \longrightarrow\A'$
for which there exists no isomorphism $\bar{i} : \B \longrightarrow \B'$ such that $\bar{i}
\circ e_A =e_{A'} \circ i$.
\item There is an $\A\in \RCA_{\alpha}$ that does not have a universal map with respect to the functor $\Nr_{\alpha}:\L\to \RCA_{\alpha}$
\end{enumroman}
\end{theorem}
\begin{demo}{Proof}(i) is proved in \cite{n}. The idea of the proof is that if (i) did not happen then $\RCA_{\alpha}$ would have the amalgamation property
which is not the case as proved by Pigozzi. (ii) follows from (i) by noting that 
\end{demo}
\begin{displaymath}
\xymatrix{
\mathfrak{B}  &\ar[l]_{e_{A}} \mathfrak{A} \ar[d]^{i} \\
              \mathfrak{B}'  & \ar[l]^{e_{A'}}\mathfrak{A}'}
\end{displaymath}
$e_A:\A\to \Nr_{\alpha}\B$ and $e_A'\circ i:\A\to \Nr_{\alpha}\B'$ for which there does not exist an isomorphism 
$\bar{f}:\B\to \B'$, such that $f\circ e_A=e_A'\circ i$.

This means that, dually, the same algebra can generate non - isomorphic algebras in extra dimensions.


We now show that this categorical result is intimately connected to various 
amalgamation properties in various classes of representable algebras.

We start with a definition:

\begin{definition}\begin{enumarab}
\item  Let $K$ be a class of algebras having a boolean reduct. 
$\A_0\in K$ is in the amalgamation base of $K$ if for all $\A_1, \A_2\in K$ and monomorphisms $i_1:\A_0\to \A_1,$ $i_2:\A_0\to \A_2$ 
there exist $\D\in K$
and monomorphisms $m_1:\A_1\to \D$ and $m_2:\A_2\to \D$ such that $m_1\circ i_1=m_2\circ i_2$. 
If in addition, $(\forall x\in A_j)(\forall y\in A_k)
(m_j(x)\leq m_k(y)\implies (\exists z\in A_0)(x\leq i_j(z)\land i_k(z) \leq y))$
where $\{j,k\}=\{1,2\}$, then we say that $\A_0$ lies in the super amalgamation base of $K$. Here $\leq$ is the boolean order.
\item $K$ has the (super) amalgamation property $((SUP)AP)$, if the (super) amalgamation base of $K$ coincides with $K$.
\end{enumarab}
\end{definition}

One can find such examples satisfying $(i)$ in theorem \ref{amal} in algebras 
that cannot be amalgamated over a common subalgebra. In fact the common subalgebra can be shown to be
the required example.
\newpage
\begin{example}
This is a family of examples. Let $i:\A_0\to \A_1, j: \A_0\to \A_2$ be monomorphisms that do not amalgamate.
By a result of Pigozzi such algebras exists among the representable algebras. 
Then $\A_0$ satisfes (i) above. If not then we can find an amalgam as follows:

\bigskip
\bigskip
\bigskip
\bigskip
\bigskip
\bigskip
\bigskip
\bigskip
\bigskip
\begin{picture}(20,0)(-70,70)

\thicklines \put(10,10){\vector(1,1){40}}\put (-5,0){$A_1$} \put
(-5,30){$k\circ e_1$}

\put (60,50){$Nr_\alpha D^+$}\put (80,70){\vector(0,1){80}}\put
(60,105){$Id$}\put (75,160){$D^+$}

\put(140,10){\vector(-1,1){40}} \put (130,30){$h\circ e_2$} \put (150,0){$A_2$}

 \put(170,0){\vector(1,0){80}} \put (200,5){$e_{2}$} \put (260,0){$ A_2^+$}


\put (145,-95){$j(A_0)$} \put(170,-90){\vector(1,0){80}} \put
(190,-100){$e_{2}\upharpoonright j(A_0)$}
\put(150,-75){\vector(0,1){60}} \put(135,-60){$Id$}

 \put (260,-100){$ Sg^{A_2^+}
(e_{2}j(A_0))$}
\put(265,-75){\vector(0,1){60}} \put(245,-60){$Id$}

\put (80,-150){$A_0$} \put(80,-160){\vector(0,-1){70}} \put
(65,-190){$e_0$} \put (80,-260){$ A_0^+$}

\put(100,-250){\vector(1,1){140}}\put (170,-190){$\bar{j}$}

\put (120,-130){$j$} \put(100,-140){\vector(1,1){40}}

\put(60,-250){\vector(-1,1){140}}\put (-15,-200){$\bar{i}$}

 \put(35,-130){$i$} \put(60,-140){\vector(-1,1){40}}

 \put(-10,0){\vector(-1,0){80}} \put (-50,5){$e_{1}$}
 \put (-120,0){$ A_1^+$}


 \put(-10,-90){\vector(-1,0){80}}\put (-5,-95){$i(A_0)$}
 \put(-60,-100){$e_{1}\upharpoonright i(A_0)$}

\put(0,-75){\vector(0,1){60}} \put(-20,-60){$Id$}


  \put (-160,-100){$Sg^{A_1^+} (e_{1}i(A_0))$}

\put(-110,-75){\vector(0,1){60}} \put(-135,-60){$Id$}

\put(-100,15){\vector(1,1){140}}\put (-50,80){$k$}

\put(270,15){\vector(-1,1){140}}\put (210,80){$h$}

\put (80,-280) {\makebox (0,0){{%
\large{  Figure 1}}}}

\end{picture}
\newpage
\end{example}

Admittedly the diagram is complicated but the idea is simple. $\A_0$ embeds into $i(\A_0)\subseteq \A_1$. 
The isomorphism can be lified to $\bar{i}$. Similarly $j$ can be lifted to $\bar{j}.$ 
We find an amlagam in the big diagram (since $\bold L$ has $AP$), and we return to the original one using the 
neat reduct functor. The property expressed in theorem \ref{amal}, {\it not holding}, 
is used to show that the isomorphisms $i$ and $j$ lift to $\bar{i}$ and $\bar{j}$.

In more detail, let $(*)$ abbreviate the negation of theorem \ref{amal} (i).
$i:\A_0\to \A_1$ is an embedding, which factors through $i:\A_0\to i(\A_0)$ and $Id:i(\A_0)\to \A_1$ and so is $j:\A_0\to \A_2$
which factors through $\A_0\to i(\A_0)$ and $Id: i(\A_0):\to \A_2$. $\A_1$ neatly embeds in $\A_1^+$ via $e_1:\A\to \Nr_{\alpha}\A^+$. 
$i(A_0)$ is a generating subreduct of 
$\Sg^{\A^+}(e_1(i(\A_0))).$  By (*) $i$ lifts to $\bar{i}:\A_0^+\to \Sg^{\A_1^+}(e_1(i(\A_0))).$ Similarly $j$ lifts to $\bar{j}:\A_0^+\to \Sg^{\A_2^+}(e_2(i(\A_0))).$ 
Now look at the big diagram. We can assume that $\A_0^+, \A_1^+$ and $\A_2^+$ are in $\bold L$; 
indeed no generality is lost if we assume that $\A_0^+=\Sg^{\A_0^+}\A_0$, $\A_1^+=\Sg^{\A_1^+}\A_1$, 
and $\A_2^+=\Sg^{\A_2^+}\A_2$. Now $\D^+$ is an amalgam of the outer diagram via $h,k$ 
and so $\Nr_{\alpha}\D^+$ is an amalgam of the inner diagram via $k\circ e_1$ and $h\circ e_2$.

\begin{corollary} Let $\bold L=\{\A\in \RCA_{\alpha+\omega}: \A=\Sg^{\A}\Nr_{\alpha}\A\}$. Then the neat reduct functor $\Nr_{\alpha}$ 
from $\bold L$ to $\RCA_{\alpha}$ with morphisms restricted to monomorphisms
does not have a right adjoint.
\end{corollary}
\begin{corollary} If $\A$ has a universal map with respect to the above functor, then $\A$ belongs to the amalgamation base of $\RK_{\alpha}$
\end{corollary}
For $\A\in \PA_{\alpha}$ a polyadic algebra and $\beta>\alpha$, a $\beta$ dilation of $\A$ is an algebra $\B\in \PA_{\beta}$
such that $\A\subseteq \Nr_{\alpha}\B.$ $\B$ is a minimal dilation of $\A$ if $A$ generates $\B.$
Let $\bold L=\{\A\in \PA_{\beta}: \Sg\Nr_{\alpha}A=A\}$. Then $\Nr_{\alpha}:\bold L\to \PA_{\alpha}$ is an equivalence.
To prove this we first note that polyadic algebras do not satisfy $(i)$ of \ref{amal}. 
But before that we need a lemma. For $X\subseteq A$, 
$\Ig^{A}X$ denotes the ideal generated by $A$.:
\begin{lemma} Let $\alpha<\beta$ be infinite ordinals. Let $\B\in \PA_{\beta}$ and $\A\subseteq \Nr_{\alpha}\B$.
\begin{enumarab}
\item if $A$ generates $\B$ then $\A=\Nr_{\alpha}\B$
\item If $A$ generates $\B$, and $I$ is an ideal of $\B$, then $\Ig^{\B}(I\cap A)=I$
\end{enumarab}
\end{lemma}\ref{net}
\begin{demo}{Proof}
\begin{enumarab}
\item Let $\A\subseteq \Nr_{\alpha}\B$ and $A$ generates $\B$ then $\B$ consists of all elements ${\sf s}_{\sigma}^{\B}x$ such that 
$x\in A$ and $\sigma$ is a transformation on $\beta$ such that
$\sigma\upharpoonright \alpha$ is one to one \cite{DM} theorem 3.3 and 4.3.
Now suppose $x\in \Nr_{\alpha}\Sg^{\B}X$ and $\Delta x\subseteq
\alpha$. There exists $y\in \Sg^{\A}X$ and a transformation $\sigma$
of $\beta$ such that $\sigma\upharpoonright \alpha$ is one to one
and $x={\sf s}_{\sigma}^{\B}.$  
Let $\tau$ be a 
transformation of $\beta$ such that $\tau\upharpoonright  \alpha=Id
\text { and } (\tau\circ \sigma) \alpha\subseteq \alpha.$ Then
$x={\sf s}_{\tau}^{\B}x={\sf s}_{\tau}^{\B}{\sf s}_{\sigma}y=
{\sf s}_{\tau\circ \sigma}^{\B}y={\sf s}_{\tau\circ
\sigma\upharpoonright \alpha}^{\A'}y.$
\item Let $x\in \Ig^{\B}(I\cap A)$. Then ${\sf c}_{(\Delta x\sim \alpha)}x\in \Nr_{\alpha}\B=\A$, 
hence in $I\cap A$. But $x\leq {\sf c}_{(\Delta x\sim \alpha)}x$, and we are done.
\end{enumarab}
\end{demo}
The previous lemma fails for cylindric algebras in general \cite{n}, but it does hold for $\Dc_{\alpha}$'s, see theorem 2.6.67, and 
2.6.71 in \cite{HMT1}.

\begin{theorem} Let $\alpha<\beta$ be infinite ordinals. Assume that $\A,\A'\in \PA_{\alpha}$ and $\B,\B'\in \PA_{\beta}.$
If $\A\subseteq \Nr_{\alpha}\B$ and $\A\subseteq \Nr_{\alpha}\B'$ and $A$ generates both then $\B$ and $\B'$ are isomorphic, 
then $\B$ and $\B'$ are isomorphic with an isomorphism that fixes $\A$ pointwise.
\end{theorem}

\begin{demo}{Proof} \cite{HMT2} theorem 2.6.72. 
We prove something stronger, we assume that $\A$ embeds into $\Nr_{\alpha}\B$ and similarly for $\A'$.
So let  $\A, \A' \in \PA_{\alpha}$ and $\beta>\alpha$. Let 
$\B, \B' \in \PA_{\beta}$ and assume that $e_A, e_{A'}$ are embeddings from $\A, \A'$ into   $\Nr_\alpha \B,
\Nr_\alpha \B'$, respectively, such that
$ \Sg^\B (e_A(A)) = \B$
and $ \Sg^{\B'} (e_{A'}(A')) = \B',$
and let $ i : \A \longrightarrow \A'$ be an isomorphism.
We need to ``lift" $i$ to $\beta$ dimensions.
Let $\mu=|A|$. Let $x$ be a  bijection  from $\mu$ onto $A.$ 
Let $y$ be a bijection from $\mu$ onto $A'$,
such that $ i(x_j) = y_j$ for all $j < \mu$.
Let $\D = \Fr_{\mu} \PA_{\beta}$ with generators $(\xi_i: i<\mu)$. Let $\C = \Sg^{\Rd_\alpha \D} \{ \xi_i : i < \mu \}.$
Then $\C \subseteq \Nr_\alpha \D,\  C \textrm{ generates } \D
~~\textrm{and so by the previous lemma }~~ \C =\Nr_{\alpha}\D.$
There exist $ f \in Hom (\D, \B)$ and $f' \in Hom (\D,
\B')$ such that
$f (g_\xi) = e_A(x_\xi)$ and $f' (g_\xi) = e_{A'}(y_\xi)$ for all $\xi < \mu.$
Note that $f$ and $f'$ are both onto. We now have
$e_A \circ i^{-1} \circ e_{A'}^{-1} \circ ( f'\upharpoonleft \C) = f \upharpoonleft \C.$
Therefore $ Ker f' \cap \C = Ker f \cap \C.$
Hence by $\Ig(Ker f' \cap \C) = \Ig(Ker f \cap \C).$
So, again by the the previous lemma, $Ker f'  = Ker f.$
Let $y \in B$, then there exists $x \in D$ such that $y = f(x)$. Define $ \hat{i} (y) = f' (x).$
The map is well defined and is as required.

\bigskip

\begin{displaymath}
    \xymatrix{
       D \ar[r]^f \ar[dr]_{f'} & B \ar[d]^{\hat{i}}  &\ar[l]_{e_{A}} A \ar[d]^{i} \\
             & B'  & \ar[l]^{e_{A'}}\mathcal{A}'}
\end{displaymath}

\end{demo}
\begin{corollary}\label{net} Let $\A, \A', i, e_A, e_{A'},$ $\B$ and $\B'$ be as in the previous proof. 
Then if $i$ is a monomorphism form $\A$ to $\A'$, then it lifts to a monomorphism $\bar{i}$ from $\B$ to $\B'$.
\end{corollary}
\begin{displaymath}
\xymatrix{
B \ar[d]^{\hat{i}}  &\ar[l]_{e_{A}} A \ar[d]^{i} \\
              B'  & \ar[l]^{e_{A'}}\mathcal{A}'}
\end{displaymath}

\begin{demo}{Proof} Consider $i:\A\to i(\A)$. Take $\C=\Sg^{\B'}(e_{A'}i(A))$. Then $i$ lifts to an isomorphism $\bar{i}\to \C\subseteq \B$.
\end{demo}
\begin{theorem} Let $\beta>\alpha$. Let $\bold L=\{\A\in \PA_{\beta}: \A=\Sg^{\A}\Nr_{\alpha}\A$\}. Let $\Nr:\bold L\to \PA_{\alpha}$ be the neat reduct functor. 
Then $\Nr$ is invertible. That is, there is a functor $G:\PA_{\alpha}\to \bold L$ and natural isomorphisms
$\mu:1_{\bold L}\to G\circ \Nr$ and $\epsilon: \Nr\circ G\to 1_{\PA_{\alpha}}$.
\end{theorem}
\begin{demo}{Proof} The idea is that a full, faithful, dense functor is invertible, \cite{cat} theorem 1.4.11. 
Let $L$ be a system of representatives for isomorphism on $Ob(\bold L)$.
For each $\B\in Ob(\PA_{\alpha})$ there is a unique $\G(B)$ in $L$ such that $\Nr(G(\B))\cong \B$.
$G(\B)$ is a minmal dilation of $\B$. Then $G:Ob(\PA_{\alpha})\to Ob(\bold L)$ is well defined. 
Choose one isomorphism $\epsilon_B: \Nr(G(B))\to \B$. If $g:\B\to \B'$ is a $\PA_{\alpha}$  morphism, then the square

\begin{displaymath}
    \xymatrix{ \Nr(G(B)) \ar[r]^{\epsilon_B}\ar[d]_{\epsilon_B^{-1}\circ g\circ \epsilon_{B'}} & B \ar[d]^g \\
               \Nr(G(B'))\ar[r]_{\epsilon_{B'}} & B' }
\end{displaymath}
commutes. By corollary \ref{net}, there is a unique morphism $f:G(\B)\to G(\B')$ such that $\Nr(f)=\epsilon_{\B}^{-1}\circ g\circ \epsilon$.
We let $G(g)=f$. Then it is easy to see that $G$ defines a functor. Also, by definition $\epsilon=(\epsilon_{\B})$ 
is a natural isomorphism from $\Nr\circ G$ to $1_{\PA_{\alpha}}$.
To find a natural isomorphism from $1_{\bold L}$ to $G\circ \Nr,$ observe that $e_{FA}:\Nr\circ G\circ \Nr(\A)\to \Nr(\A)$ is an isomorphism.
Then there is a unique $\mu_A:\A\to G\circ \Nr(\A)$ such that $\Nr(\mu_{\A})=e_{FA}^{-1}.$
Since $\epsilon^{-1}$ is natural for any $f:\A\to \A'$ the square
\bigskip
\bigskip
\begin{displaymath}
    \xymatrix{ \Nr(A) \ar[r]^{\epsilon_{\Nr(A)}^{-1}=\Nr(\mu_A)}\ar[d]_{\Nr(f)} & \Nr\circ G\circ \Nr(A) \ar[d]^{\Nr\circ G\circ \Nr(f)} \\
               \Nr(A')\ar[r]_{\epsilon_{FA}^{-1}=\Nr(\mu_{A'})} & \Nr\circ G\circ \Nr(A') }
\end{displaymath}

commutes, hence the square

\bigskip
\begin{displaymath}
    \xymatrix{ A \ar[r]^{\mu_A}\ar[d]_f & G\circ \Nr(A) \ar[d]^{G\circ \Nr(f)} \\
               A'\ar[r]_{\mu_{A'}} & G\circ \Nr(A') }
\end{displaymath}

commutes, too. Therefore $\mu=(\mu_A)$ is as required.
\end{demo}

Let $C$ be the reflective subcatogory of $\RCA_{\alpha}$ that has universal maps. Then $\Dc_{\alpha}\subseteq \bold L$.
And indeed we have:
\begin{theorem}\label{SUP} Let $\alpha\geq \omega$ . Let $\A_0\in \Dc_{\alpha}$, $\A_1,\A_2\in \RCA_{\alpha}$
and $f:\A_0\to \A_1$ and $g:\A_0\to \A_2$ be monomorphisms. Then there exists $\D\in \Nr_{\alpha}\CA_{\alpha+\omega}$
and $m:\A_1\to D$ and $n:\A_2\to \D$ such that $m\circ f=n\circ g$. Furthermore $\D$ is a super amalgam.
\end{theorem}
\begin{demo}{Proof}
$\Dc_{\alpha}$ does not satify property (i) in theorem \ref{amal}. 
Looking at figure 1, assuming that the base algebra $\A_0$ is in $\Dc_{\alpha}$, we obtain 
$\D\in \Nr_{\alpha}\CA_{\alpha+\omega}$
$m:\A_1\to \D$, and $n:\A_2\to \D$
such that $m\circ i=n\circ j$.
Here $m=k\circ e_1$ and $n=h\circ e_2$.
Denote $k$ by $m^+$ and $h$ by $n^+$.
Now we further want to show that if $m(a) \leq n(b)$,
for $a\in A_1$ and $b\in A_2$, then there exists $t \in A_0$
such that $ a \leq i(t)$ and $j(t) \leq b$.
So let $a$ and $b$ be as indicated . We have  $m^+ \circ e_1(a) \leq n^+ \circ e_2(b),$ so
$m^+ ( e_1(a)) \leq n^+ ( e_2(b)).$
Since $\bold L$ has $SUPAP$, there exist $ z \in A_0^+$ such that $e_1(a) \leq \bar{i}(z)$ and
$\bar{j}(z) \leq e_2b)$.
Let $\Gamma = \Delta z \smallsetminus \alpha$ and $z' =
c_{(\Gamma)}z$. (Note that $\Gamma$ is finite.) So, we obtain that
$e_1(c_{(\Gamma)}a) \leq \bar{i}(c_{(\Gamma)}z)~~ \textrm{and} ~~ \bar{j}(c_{(\Gamma)}z) \leq
e_2(c_{(\Gamma)}b).$ It follows that $e_A(a) \leq \bar{i}(z')~~
\textrm{and} ~~ \bar{j}(z') \leq e_B(b).$ Now $z' \in \Nr_\alpha \A_0^+
= \Sg^{\Nr_\alpha \A_0^+} (e_{A_0}(A_0)) = A_0.$ Here we use \cite{HMT1} 2.6.67.
So, there exists $t \in C$ with $ z' = e_C(t)$. Then we get
$e_1(a) \leq \bar{i}(e_0(t))$ and $\bar{j}(e_1(t)) \leq e_2(b).$ It follows that $e_1(a) \leq e_A \circ i(t)$ and
$e_2 \circ j(t) \leq
e_2(b).$ Hence, $ a \leq i(t)$ and $j(t) \leq b.$
We are done.
\end{demo}
\begin{corollary} For $\alpha\geq \omega$, $\Dc_{\alpha}$ is contained in the superamalgamation base of $\RCA_{\alpha}$
\end{corollary}
\begin{corollary} If $\A$ has universal maps and satisfies $NS$ then $\A$ lies in the $SUPAP$ base
\end{corollary}

Call a system of varieties neat if it is a system of varieties definable by schemes
satisfying the the finiteness generating condition, 
and satisfying that for all $\A\in K_{\alpha}$ there exists 
$\B\in K_{\alpha+\omega}$ such that for all $X\subseteq A$, $\Sg^{\A}X=\Nr_{\alpha}\Sg^{\B}X$. 

Call a system of varieties nice  if the neat reduct functor has a right adjoint, and 
$K_{\alpha}=Kn_{\alpha}$ have $SUPAP$, and each $K_{\omega}$ is axiomatized by a finite schema. 
Is there a neat or /and nice system of varieties definable 
by (finitely many) schemas? 
This is a difficult question that lies at the heart of the process of algebriasation, and is strongly related to the so called finitizability 
problem in algebraic logic. It basically asks whether the is an optimal combination of the cylindric paradigm and the polyadic one; optimal here meaning
that it avoids negative properties of both and, on the other hand, shares their positive properties. If we do not insist on 'definable by schemes' 
then there is such a system
\cite{Sain}.

\section*{Some Remarks}

\begin{enumarab}
\item  For cylindric algebras a single algebra can generate two non isomorphic algebras in extra dimensions
and dually two non isomporhic algebras can generate the same algebra in extra dimensions. This cannot happen for polyadic algebras.
For representable  cylindric algebras this can be formulated as follows. If one takes a minimal dilation of an algebra
and then apply the neat reduct functor, one does not  necessarily end where he started.
This again cannot happen for the polyadic case.

\item In the introduction of \cite{HMT1} it was asked by Henkin and Monk whether theorems 6.67-2.6.71, 2.6.72 which hold in $\Dc_{\alpha}$
continue to hold for arbitary $\RCA_{\alpha}$. Category theory was then not mature enough. Now it can be paraphrased using categorical jargon.
Their queston, as shown herein,  can be reformulated as to whether the neat reduct functor has a right adjoint.

\item It is known from several results in the literature that cylindric algebras and polyadic algebras belong to different paradigms.
For example the class of representable cylindric algebras cannot be axiomatized
by finitely many schema and fails to have $AP$, while in contrast polyadic algebras are all representable, and the class has
$SUPAP$. The possesing of a right adjoint (in fact an equivalence) for the neat reduct functor for $\PA$'s and its absence in representable 
$\CA$'s highlights another difference which is crucial in proving amalgamation results.

\item The interaction between the theories of polyadic algebras and cylindric algebras have been recently studied with pleasing progress 
and a plathora of results, to mention a few references in this connection \cite{ANS}, \cite{F1}, \cite{F2}, \cite{FM}, \cite{trans}, \cite{Fer4}
and \cite{Sagi}.

\item $\Dc_{\alpha}$ and $\PA_{\alpha}$ have a lot in common. In fact, $\Dc_{\alpha}$ is contained in the $SUPAP$ base of $\RCA_{\alpha}$ by theorem 
\ref{SUP} and $ \PA_{\alpha}$ has $SUPAP$. The proof of the latter is much more involved \cite{SUP}. 
The neat reduct functor restricted to both have a right adjoint. So $\Dc_{\alpha}$, it seems, is the largest subclass of $\RCA_{\alpha}$ 
which enjoys positive properties of 
$\PA_{\alpha}$
\item Let $\FPA_{\alpha}$ be the class of algebras that are like cylindric algebras in that cylindrifiers are only finite and we have all substitutions.
This is a variety that can be proved using the techniques in \cite{SUP} to have $SUPAP$.
$\FPA_{\alpha}$ has a double behaviour. In one of its facets it resembles cylindric algebras with only finite cylindrifiers available, 
and on the other hand it has all substitutions available 
a property it shares with the polyadic paradigm. $\FPA_{\alpha}$ can be formulated as a system of varieties that is both neat and nice.

\item Directed cylindric algebras introduced by N\'emeti \cite{Sag} belong to the polyadic paradigm! Indeed it can be formulated 
as a system of varieties (of finite dimensions) that are neat and nice \cite{Sag}, \cite{Andras}.
So are the algebras investigated by Sain in \cite{Sain}. Both can be considered as solutions to the finitizability problem, the former in finite dimensions, 
the latter in infinite dimensions. For those neat reducts commute with forming subalgebras.

\item It is no coincidence that $SUPAP$ and the invertibility of the neat reduct functor come together. For the invertibility of the neat reduct functor, 
roughly says, that terms definable in infinitely many extra dimensions are alraedy term definable.
This is a form of definability, which as we have shown here is closely linked to classical definability results, 
like Craig interpolation and Beth definability.

\end{enumarab}

\end{document}